\documentstyle[epsfig,11pt]{article}
\input epsf
\onecolumn 
\begin{document}

%\input{symbols}%%%%%%%%%%%%%%%%%%%%%%%%%%%%%%%%%%%%%%%%%%%%%%%%%%%%%%%%%%%
%for the meaning of the symbols see notation.tex
\newcommand {\Pind} {\hspace{.51cm}}  % paragraph indent 
\newcommand {\vare}{\varepsilon}
\newcommand {\grad}{\bigtriangledown}
\newcommand {\RR}{\,\hbox{\rm R}\!\!\!\!\!\hbox{\rm I}\,\,\,}
\newcommand {\NN}{\,\hbox{\rm N}\!\!\!\!\!\hbox{\rm I}\,\,\,}
\newcommand {\implies}{Longrightarrow}
\newcommand {\ProofEnd} {\hfill \nobreak \square \medbreak}
\newcommand{\square}{\hbox{${\vcenter{\hrule height.4pt
  \hbox{\vrule width.4pt height6pt \kern6pt
     \vrule width.4pt}
  \hrule height.4pt}}$}}

\newtheorem {theorem} {Theorem}
\newtheorem {proposition} [theorem] {Proposition}
\newtheorem {lemma} [theorem] {Lemma}
\newtheorem {definition} [theorem] {Definition}
\newtheorem {corollary} [theorem] {Corollary}
\newtheorem {remark} [theorem] {Remark}
\newtheorem {note} [theorem] {Note}
\newtheorem {example} [theorem] {Example} 

%%%%%%%%%%%%%%%%%%%%%%%%%%%%%%%%%%%%%%%%%%%%%%%%%%%%%%%%%%%%%%%%%%%%%%%%%%%%

\begin{titlepage} 
\begin{center}
{\huge\bf Nonnegative Ricci curvature, small linear diameter growth
 and finite generation of fundamental groups\\}

\vspace*{7ex} 
{\Large\bf {C. Sormani, 1998} }

\vspace{1cm}
{\large \bf {sormani@math.jhu.edu}}
\end{center}

\end{titlepage} 
%%%%%%%%%%%%%%%%%%%%%%%%%%%%%%%%%%%%%%%%%%%%%%%%%%%%%%%%%%%%%%%%%%%%%%%%%%%%
\pagestyle{plain}  %try

\pagenumbering{arabic}

\begin{section} {Introduction}

\Pind
In 1968, Milnor conjectured that a complete noncompact manifold, $M^n$, with
nonnegative Ricci curvature has a finitely generated fundamental
group [Mi].  This was proven for a manifold with nonnegative sectional
curvature by Cheeger and Gromoll [ChGl1].  However, it remains an
open problem even for manifolds with strictly positive Ricci curvature.  

The conjecture is of particular interest because, if it is true, then by 
work of Cheeger-Gromoll, Milnor and Gromov, the fundamental group is almost 
nilpotent [ChGl2] [Mi] [Gr].   
On the other hand, given any finitely generated torsion free nilpotent group,
Wei has constructed an example of a manifold with positive Ricci
curvature that has the given group as a fundamental group [Wei]. 

Schoen and Yau have proven the conjecture in dimension 3 for 
manifolds with strictly positive Ricci curvature [SchYau].  In fact they have
proven that such a manifold is diffeomorphic to Euclidean space.

Anderson and Li have each proven that if a manifold
with nonnegative Ricci curvature has Euclidean volume growth, then the
fundamental group is actually finite [And] [Li].  Anderson uses volume
comparison arguements while Li uses the heat equation to prove this theorem.

Abresch and Gromoll have proven that manifolds with 
small diameter growth, $(o(r^{1/n}))$, nonnegative Ricci curvature and 
sectional curvature bounded away from negative infinity have finite 
topological type [AbrGrl, Thm A].  Thus the fundamental group is finitely 
generated in such manifolds.    Their theorem is proven using an inequality
refered to as the Excess Theorem, [AbrGrl, Prop. 2.3], which is
a crucial ingredient in this paper as well.

Abresch and Gromoll have also proven that a manifold 
with nonnegative Ricci curvature has linear almost intrinsic diameter growth 
[AbrGrl, Prop 1.1].  That is, given any $\vare >0$, there exists an explicit 
constant, $C(\vare, n)$, such that given any $p$ and $q$ in
$\partial B_{x_0}(r)$ joined by a curve in the annulus, 
$Ann_{x_0}(r-\vare r, r+ \vare r)$, if 
$\sigma$ is the shortest such curve then $L(c) \le C(\vare, n) r$.  
This proposition is proven using volume comparison arguements.

In this paper we prove that a manifold with small linear diameter growth
has a finitely generated fundamental group.     

\begin{theorem} \label{diamgroup}
There exists a universal constant,
\begin{equation}\label{eqnSn}
 S_n = \frac{n}{n-1}\frac 1 {4} \frac 1 {3^n} 
\left(\frac {n-2}{n-1}\right)^{n-1},
\end{equation}
such that if $M^n$ is complete and noncompact with nonnegative Ricci curvature
and has  small linear diameter growth, 
$$
\limsup_{r\to\infty} \frac{diam(\partial B_p(r))}{r} < 4 S_n,
$$
then it has a finitely generated fundamental group.
\end{theorem}

Note that unlike the theorem of Abresch and Gromoll we only control
the fundamental group.   However, we do not require a lower sectional curvature
bound and the requisite diameter control is much weaker. 

We also prove a more general theorem, 
Theorem~\ref{PoleGroup}.  If $M^n$ has nonnegative Ricci 
curvature and an infinitely generated fundamental group then it has a tangent 
cone at infinity which is not polar.  Currently it is not known
whether a manifold with nonnegative Ricci curvature can have a tangent
cone at infinity which is not polar [ChCo].

The definitions of tangent cone and polar length spaces are reviewed
in Section 5 [Defn~\ref{DefnTanCone}, 
Defn~\ref{DefnPole}] along with the precise statement of 
Theorem~\ref{PoleGroup}.  The precise statement involves
the universal constant, $S_n$, defined in (\ref{eqnSn}). 

Note that if $M^n$ has Euclidean volume growth then by Cheeger and
Colding, [ChCo1], its tangent cones at infinity have poles, and thus $M^n$
has a finitely generated fundamental group.  However, in this case, 
Anderson and Li have each proven that the fundamental group is actually 
finite.  [And, Cor 1.5], [Li].

On the other hand, if $M^n$ has linear volume growth then by [So],
its diameter growth is sublinear and its tangent cone at infinity is 
$[0,\infty)$ or $(-\infty, \infty)$, thus we have the following
corollary of either theorem.

\begin{corollary} \label{volgroup}
If $M^n$ is complete with nonnegative Ricci curvature and
linear volume growth then it
has a finitely generated fundamental group.
\end{corollary}

The proofs of both Theorem~\ref{diamgroup} and Theorem~\ref{PoleGroup} 
are based on two main lemmas.  

The Halfway Lemma, [Lemma~\ref{Halfway}], concerns complete Riemannian
manifolds with infinitely generated 
fundamental groups.  It does not require the Ricci curvature condition and 
makes a special selection of generators and representative loops, such that
the loops are minimal halfway around.  It is stated and proven in Section 2.

The Uniform Cut Lemma, [Lemma~\ref{unifcut}], gives a uniform 
estimate on special cut points which are the halfway points of noncontractible
geodesic loops in a manifold with nonnegative Ricci curvature.  The proof 
applies the Excess Theorem of Abresch and Gromoll.  It appears in Section 3.

Section 4 contains the proof of Theorem~\ref{diamgroup} and Section 5
contains the proof and background material for Theorem~\ref{PoleGroup}.
Both of these sections require the results of Sections 1 and 2
but are independant of each other.

The author would like to thank Professor Cheeger for suggesting that
Theorem~\ref{diamgroup} be strenthened to its current form and for
his assistance during the revision process. 
Background material can be found in [Ch] and [GrLaPa].

\end{section}

\begin{section} {The Halfway Lemma}

In this section our manifold, $N^n$, is a complete Riemannian manifold 
but does not have a bound on its Ricci curvature.  It may or may not
be noncompact.   Recall that a group
is infinitely generated if it does not have a presentation with a finite
set of generators.  If $\pi_1(N)$ is not finitely
generated, $N^n$ must be noncompact.

\begin{definition} {\em
If $G$ is a group, we say that $\{g_1, g_2,...\}$
is a ordered set of {\em independant generators} of $G$ if each $g_i$ can not
be expressed as a word in the previous generators and their inverses,
$g_1, g_1^{-1}, ... g_{i-1}, g_{i-1}^{-1}$.}
\end{definition}

\begin{definition} \label{MinRep} \label{MinRepLength}
{\em Given $g\in \pi_1(N)$, we say $\gamma$ is a {\em minimal representative
geodesic loop} of
$g$ if $\gamma = \pi \circ \tilde{\gamma}$,
where $\tilde{\gamma}$ is minimal from $\tilde{x_0}$ to $g \tilde{x_0}$
in $\tilde{N}$.  Note that 
$L(\gamma)=d_{\tilde{M}}(\tilde{x_0}, g\tilde{x_0})$.}
\end{definition}

We now state the Halfway Lemma.

\begin{lemma} \label{Halfway} [Halfway Lemma]  %(proven 7/1/97 p 3-4)
Let $x_0 \in N^n$ where $N^n$ is a complete Riemannian manifold 
with a fundamental group $\pi_1(N,x_0)$. 
Then there exists an ordered set of independant generators 
$\{g_1, g_2, g_3...\}$ of  $\pi_1(N, x_0)$ with minimal representative 
geodesic loops, $\gamma_k$, of length $d_k$ such that
\begin{equation}
d_N(\gamma_k(0), \gamma_k(d_k/2))= d_k/2.
\end{equation}
If $\pi_1(N)$ is infinitely generated we have a sequence of such generators.
\end{lemma}

%Warning: cannot just take a subgroup of the fund group or trouble at * below.

\begin{definition} \label{HalfGen} {\em
Given $x_0 \in N$, we will call an ordered set of generators 
of $\pi_1(N, x_0)$ as constructed in Lemma~\ref{Halfway}, a set of
{\em halfway generators} based at $x_0$ of $\pi_1(N)$. }
\end{definition}

\vspace{.5cm}
\noindent {\bf Proof:}
In order to prove the Halfway Lemma, we need to choose a sequence of
generators whose representative curves don't have redundant extra looping.
Let $G=\pi_1(N^n, x_0)$.

%INSERT DIAGRAM 1

We first define the
sequence of generators $\{g_1, g_2, g_3...\}$.
Define $g_1\in G$ such that it minimizes distance,
$$
d_{\tilde{N}}(\tilde{x_0}, g_1\tilde{x_0})\le
d_{\tilde{N}}(\tilde{x_0}, g\tilde{x_0}) \qquad \forall g\in G.
$$
Let $G_i \in G$ be the elements of G generated by 
$g_1, ... g_i$ and their inverses.  So $G_1=\{e,g_1, g_1^{-1}, g_1^2,...\}$.
Define each $g_k \in G$ iteratively such that each  minimizes distance
among all elements in $G \setminus G_{k-1}$, 
$$
d_{\tilde{N}}(\tilde{x_0}, g_k\tilde{x_0})\le
d_{\tilde{N}}(\tilde{x_0}, g\tilde{x_0}) \qquad \forall g\in 
G\setminus G_{k-1}.
$$
Note that $G\setminus G_{k-1}$ is nonempty for all $k$ when we have an 
infinitely generated fundamental group.  

Thus if there exists $h \in G$ with
\begin{equation} \label{noshort}
d_{\tilde{N}}(\tilde{x_0}, h\tilde{x_0}) < 
d_{\tilde{N}}(\tilde{x_0}, g_k\tilde{x_0})
\end{equation}
then $h \in G_{k-1}$. 

Let $\tilde{\gamma_k}$ be a minimal geodesic in $\tilde{N}$ from
$\tilde{x_0}$ to $g_k \tilde{x_0}$.  Let 
$\gamma_k(t)=\pi(\tilde{\gamma_k}(t))$.  Then $\gamma_k$ has no
conjugate points on $[0, d_k)$.  

Given a curve $C:[0,d] \to N$, let \newline
$C(t_1 \to t_2)$ represent the segment of $C$ running from
$t_1$ to $t_2$.  We allow $t_2<t_1$.

Suppose that there is a $k\in \NN$ such that
\begin{equation} \label{contrahalf}
d_N(\gamma_k(0), \gamma_k(d_k/2))< d_k/2.
\end{equation}
Then there exists $T < d_k/2$ 
such that $\gamma_k(T)$ is a cut point of $x_0=\gamma(0)$.  Thus there exists
a geodesic $\sigma$ from $\gamma_k(T)=\sigma(0)$ to $x_0=\sigma(T)$.
This geodesic segment cannot overlap $\gamma_k(T \to d_k)$ because
it has a shorter length, $T < d_k/2$.

Then there exists elements 
$h_1=[\sigma(0\to T) \circ \gamma_k (0 \to T)] \in \pi_1(N)$ %*
and $h_2=[\gamma_k(T \to d_k) \circ \sigma(T \to 0)]\in \pi_1(N)$.
Furthermore, since $\sigma$ meets $\gamma$ at an angle, the lift is
not a minimal geodesic, so
$$
d_{\tilde{N}}(\tilde{x_0}, h_1\tilde{x_0}) 
< L(\sigma(0\to T) \circ \gamma_k (0 \to T))
< 2T < d_k
$$
and
$$
d_{\tilde{N}}(\tilde{x_0}, h_2\tilde{x_0}) < 
L(\gamma_k(T \to d_k) \circ \sigma(T \to 0))<
T + (d_k-T) = d_k.
$$ 
Thus, by (\ref{noshort}), $h_1$ and $h_2$ are in $G_{k-1}$. 
So $g_k =h_2 \circ h_1 \in G_k$ contradicting our choice of
$g_k \in G \setminus G_{k-1}$.  So our assumption in (\ref{contrahalf}) 
is false.
\ProofEnd

\end{section}

\begin{section} {The Uniform Cut Lemma}

In this section, $M^n$  is a complete manifold with nonnegative
Ricci curvature of dimension, $n \ge 3$ .  It may or may not be compact. 
Recall that when $n=2$, Ricci curvature is just sectional curvature
and Theorems~\ref{diamgroup} and~\ref{PoleGroup} follow from [ChrGrl1].

The Uniform Cut Lemma describes special cut points which are the halfway
points, $\gamma(D/2)$, of geodesic loops, $\gamma$, from a base point, 
$x_0$.  Recall that no geodesic from $x_0$ through a cut point of 
$x_0$ is minimal after passing through that cut point.  So if 
$x \in \partial B_{x_0}(RD)$ where $R > 1/2$, then any minimal geodesic 
from $x$ to $x_0$ does not hit the cut point $\gamma(D/2)$.  Thus 
$$
d_M(x, \gamma(D/2)) > (R\,-\,1/2) D.
$$  
The Uniform Cut Lemma gives a uniform and scale invariant estimate for 
this inequality.

\begin{lemma} [Uniform Cut Lemma] \label{unifcut}
Let $M^n$ be a complete with nonnegative Ricci curvature
and dimension $n \ge 3$.  Let $\gamma$ be a noncontractible geodesic loop 
based at a point, $x_0 \in M^n$,  of length,
$L(\gamma)=D$, such that the following conditions hold:

i) If $\sigma$ based at $x_0$ is a loop homotopic to $\gamma$ then 
$L(\sigma)\ge D$

ii) The loop $\gamma$ is minimal on $[0,D/2]$ and 
is also minimal on $[D/2, D]$. 

\noindent
Then there is a universal constant $S_n$, defined in (\ref{eqnSn}),
such that if $x \in \partial B_{x_0}(RD)$ where $R \ge (1/2 +S_n) $ then
$$
d_M(x, \gamma(D/2))\ge  (R-1/2)D + 2 S_n D.
$$
\end{lemma}

Note that the minimal representative geodesic loops of the
halfway generators constructed in the Halfway Lemma satisfy
the hypothesis of the Uniform Cut Lemma.

\vspace{.5cm}
\noindent {\bf Proof}
We first prove the lemma for $R=(1/2 +S_n)$.

Assume on the contrary, that there is a point $x\in M^n$ such that
$d_M(x, x_0)=(1/2 +S_n) D$ and $d_M(x, \gamma(D/2))=H <  3 S_n D$.
Let $C:[0,H] \longmapsto M^n$ be a minimal geodesic from
$\gamma(D/2)$ to $x$. 

Let $\tilde{M}$ be the universal cover of $M$, let $\tilde{x_0}\in \tilde{M}$
be a lift of $x_0$, and let $g\in \pi_1(M,x_0)$ be the element represented
by the given loop, $\gamma$.  By the first condition on $\gamma$, its lift, 
$\tilde{\gamma}$, is a 
minimal geodesic running from $\tilde{x_0}$ to $g\tilde{x_0}$.
Thus 
$$
d_{\tilde{M}}(\tilde{x_0},g\tilde{x_0})=D.
$$

We can lift the joined curves, $C(0 \to H) \circ \gamma(0 \to D/2)$, to 
a curve in the universal cover, $\tilde{C} \circ \tilde{\gamma}$, 
which runs from $\tilde{x_0}$ through $\tilde{\gamma}(D/2)$
to a point $\tilde{x} \in \tilde{M}$.
Note that $L(\tilde{C})=L(C)=H$.
  
We can examine the triangle formed by $\tilde{x_0}$, $g\tilde{x_0}$
and $\tilde{x}$ using the Excess Theorem of Abresch and Gromoll [AbrGrl].  
For easy reference we introduce their variables $r_0$ and $r_1$.

By our assumption on $x$, 
\begin{equation}
r_0  =  d_{\tilde{M}}(\tilde{x}, \tilde{x_0}) \ge  d_M(x, x_0)
 =  (1/2 + S_n) D.
\end{equation}
Furthermore,
$$
r_1=d_{\tilde{M}}(g_k \tilde{x}, \tilde{x_0}) 
\ge d_M(x, x_0) =(1/2 + S_n) D.
$$

The excess of $\tilde{x}$ relative to $\tilde{x_0}$ and $g\tilde{x_0}$
satisfies
\begin{equation} \label{excess1}
e(\tilde{x}):=  r_0 + r_1 -d(\tilde{x_0}, g \tilde{x_0})  \ge 
  2(1/2 + S_n) D -D = 2 S_n D.
\end{equation}

On the other hand, by the Excess Theorem [AbrGrl, Prop 2.3], we can estimate
the excess from above in terms of the distance, $l$, from $\tilde{x}$
to the minimal geodesic, $\tilde{\gamma}$.
In particular for $n \ge 3$, $Ricci \ge 0$, they have proven that
\begin{equation} \label{AbrGrl1}
e(\tilde{x}) \le
2 \left(\frac {n-1}{n-2}\right)\left(\, \frac 1 2 \, C_3 \, l^n \,
\right)^{1/(n-1)}
\end{equation}
where
\begin{equation} \label{AbrGrl2}
C_3= \frac {n-1}{n}\left( \frac 1 {r_0-l} \,+\, \frac 1 {r_1 - l}\right)
\end{equation}
if $l < \min \{r_0, r_1\}$.  

We now need to estimate $l$ from above.
Suppose that the closest point on $\tilde{\gamma}$ to $\tilde{x}$
occurs at a point $\tilde{\gamma}(t_0)$.  Then
\begin{eqnarray*}
l=d_{\tilde{M}}(\gamma(t_0), \tilde{x}) & \le &
d_{\tilde{M}}(\tilde{\gamma}(D/2), \tilde{x})\\
& \le & L(C) = HD < 3S_n D.
\end{eqnarray*}
Since $S_n <1/20$ for $n \ge 3$,
$$
r_0-l  \ge  (1/2 +S_n)D - 3S_n D > D/4. 
$$
and, similarly,
$$
r_1-l > D/4.
$$
In particular, $l < \min \{r_0, r_1\}$.

Substituting this into the Abresch and Gromoll's estimate (\ref{AbrGrl1},
\ref{AbrGrl2}), we have 
\begin{eqnarray*}
e(\tilde{x_k}) &<&
2 \left(\frac {n-1}{n-2}\right)
\left(\frac 1 2 \left(\frac {n-1} {n}\right)
\left(\frac 2 {D/4}\right)
 \left(3S_n D \right)^n \right)^{1/(n-1)}  \\
& \le & 2 D \left(\frac {n-1}{n-2}\right)
\left(4 \left(\frac {n-1} {n}\right)
 \left(3S_n  \right)^n \right)^{1/(n-1)}. 
\end{eqnarray*}

Combining this with (\ref{excess1}) and cancelling $D$, we get
\begin{equation}
2S_n < 2 \left(\frac {n-1}{n-2}\right)
\left( 4 \left(\frac {n-1} {n}\right)
 \left(3S_n  \right)^n \right)^{1/(n-1)}.
\end{equation}
Cancelling $2$ and exponentiating, we have
$$
S_n^{n-1} =
4 \left(\frac {n-1}{n-2}\right)^{n-1}\frac {3^n (n-1)}{n}
S_n^n,
$$
and
$$
S_n > \frac{n}{n-1}\frac 1 {4} \frac 1 {3^n} 
\left(\frac {n-2}{n-1}\right)^{n-1}.
$$
This contradicts the definition of $S_n$ in (\ref{eqnSn}) and we've
proven the theorem for $R=(1/2 +S_n)$.

If $R \ge (1/2 +S_n)$, let $x \in \partial B_{x_0}(RD)$ and let
$y \in \partial B_{x_0}((1/2 +S_n)D)$ be a point on a minimal geodesic
from $x$ to $\gamma(D/2)$.  Then, by the above case,
\begin{eqnarray*}
d_M(x, \gamma(D/2))&=& d_M(x,y) + d_M(y, \gamma(D/2))\\
& \ge & (RD-(1/2 +S_n)D) + 3 S_n D = (R-1/2)D + 2 S_n D.
\end{eqnarray*}

\ProofEnd

\end{section}

\begin{section} {The Small Linear Diameter Growth Theorem}

\Pind
In this section we prove Theorem~\ref{diamgroup}.

\vspace{.5cm}
\noindent {\bf Proof of Theorem~\ref{diamgroup}:}
Assume that $M^n$ has an infinitely generated fundamental group
$\pi_1(M^n, x_0)$.  Then by the Halfway Lemma, [Lemma~\ref{Halfway}], 
there is a sequence of halfway generators, $g_k$, whose minimal 
representative geodesic loops based at $x_0$, $\gamma_k$, satisfy the 
hypothesis of the Uniform Cut Lemma [Lemma~\ref{unifcut}].
Let $d_k=L(\gamma_k)$.  Note that $d_k$ diverges to infinity.

By the Uniform Cut Lemma, if 
$x_k\in \partial B_{x_0}((1/2 +S_n)d_k)$  then
$$
d_M(x_k, \gamma(d_k/2))\ge  3 S_n d_k.
$$
Thus the point $y_k \in \partial B_{x_0}((1/2)d_k)$ on the minimal geodesic
from $x_k$ to $x_0$, satisfies,
\begin{eqnarray*}
d_M(y_k, \gamma_k(d_k/2)) &\ge& d_M(x_k, \gamma_k(d_k/2)) - d(x_k, y_k)\\
& \ge & (3 S_n d_k) - (S_n d_k) = 2 S_n d_k.
\end{eqnarray*}

This allows us to estimate the diameter growth,
\begin{eqnarray*}
\limsup_{r\to\infty} \frac {diam(\partial B_{x_0}(r)}{r}
& \ge & \limsup_{k\to \infty} \frac{ d(y_k, \gamma_k(d_k/2))}{(d_k/2)} \\
& \ge & \limsup_{k\to \infty} \frac{2 S_n d_k}{d_k/2}=4 S_n.
\end{eqnarray*}
This contradicts the small linear diameter growth of $M^n$.
\ProofEnd
\end{section}

\begin{section} {The Pole Group Theorem} 

\Pind
In this section we state and prove Theorem~\ref{PoleGroup}.
We begin with some background material.  

Given a complete noncompact manifold, $M^n$, with nonnegative Ricci 
curvature, we can define tangent cones at infinity by taking
any pointed sequence of rescalings of the manifold.  A subsequence
of such a sequence must converge in the Gromov Hausdorff topology
to a pointed length space, $(X, x_0)$, by Gromov's Compactness Theorem.  
[GrLaPa]  %(add appendix?)

\begin{definition} \label{DefnTanCone}  {\em
A pointed length space, $(X,x_0)$, is called a {\em tangent cone at infinity}
of $M$ if there exists a point $p\in M$ and a sequence $r_k$ of 
positive real numbers diverging to infinity such that for all $R>0$,
$$
d_{GH}\big( (B_R(x_0)\in X, x_0, d_X), (B_R(p)\in M, p, d_M/r_k) \big) \to 0,
$$
as $r_k\to\infty$.  Here $d_M$ is the length space distance function
on $M$ induced by the Riemannian metric $g_M$ and $d_{GH}$ is the 
Gromov-Hausdorff distance.}
\end{definition}

Note that $(X,x_0)$ need not be unique.  See [Per] and [ChCo, 8.37].
%FIND ORIGINAL REF FOR TAN CONE
Furthermore, $(X, x_0)$ need not be a metric cone unless the manifold
has Euclidean volume growth [ChCo].

\begin{definition} \label{DefnPole} {\em [ChCo, sect 4]
A length space, $X$, has a {\em pole} at a point $x \in X$ if
for all $y$ not equal to $x$ there exists a curve $\gamma: [0,\infty) \to X$
such that $\gamma(0)=x$, $d_X(\gamma(t), \gamma(s))=|s-t|$ for all $s,t \ge 0$,
and $\gamma(d(x,y))=y$. }
\end{definition}

There is no known example of a manifold with nonnegative Ricci curvature
with a tangent cone at infinity
which does not have a pole at its base point [ChCo].  In order to find
an example of such a manifold, intuitively one would need to construct 
a sequence of cut points on the manifold which remain {\em uniformly cut}
even after rescaling.  By Lemmas~\ref{Halfway} and~\ref{unifcut}, such cut 
points exist if the manifold has an infinitely generated fundamental group.
This is the intuition behind  
Theorem~\ref{PoleGroup} and Theorem~\ref{diamgroup}.

\begin{theorem}  \label{PoleGroup} [Pole Group Theorem] 
If a complete noncompact manifold, $M^n$, with nonnegative Ricci
curvature  has a fundamental group which is not finitely 
generated, then it has a  tangent cone at infinity, $(Y, y_0)$, which does 
not have a pole at its base point.  

In fact, if $(Z, z_0)$ is a length space with 
\begin{equation} \label{eqnZY}
d_{GH}((B_{z_0}(1), z_0, d_Z), (B_{y_0}(1), y_0, d_Y)) < S_n/4.
\end{equation}
where $S_n$ was defined in (\ref{eqnSn}), 
then $Z$ does not have a pole at $z_0$.
\end{theorem}

%Note that this theorem requires that Ricci curvature be nonnegative. 
%The 2 dimensional ladder manifold [see diagram] can be created with
%$Ricci \ge -1$,
%has an infinitely generated fundamental group and its tan cone at infinity
%is just the real line.  

\noindent {\bf Proof:}

First we choose a special sequence of rescalings of $M^n$ and a 
corresponding tangent cone at infinity.  Let $x_0$ be any base point
in $M$.  By the Halfway Lemma, there is a sequence of halfway generators, 
$g_k$, of lengths, $d_k$, corresponding to a base point $x_0$.
Take a subsequence of this sequence such that $M_k=(M^n, x_0, d_M/d_k)$
converges to a tangent cone $Y=(Y, y_0, d_Y)$.  [GrLaPa] %Better ref?

By the Halfway Lemma and the Uniform Cut Lemma, we know that for all 
$k \in \NN$, for all $r \ge 1/2 +S_n$ and for all
$x \in \partial B_{x_0}(rd_k)$  then
\begin{equation}\label{lemmaseqn}
d_M(x, \gamma_k(d_k/2))\ge  (r-1/2+2 S_n) d_k.
\end{equation}

We want to find a ``cut point'' in $Y$ and in $Z$. 
By Definition~\ref{DefnTanCone}, taking $R=1$ 
there exists $N \in \NN$ such that for all $k \ge N$ ,
$$
d_{GH}((B_{y_0}(1)\subset Y, y_0, d_Y), 
(B_{x_0}(1)\subset M_k, x_0, d_{M_k})) < S_n/12.
$$
Thus, by (\ref{eqnZY}), for all $k \ge N$, the length space $Z$ satisfies
\begin{equation} \label{eqnvareSn}
d_{GH}((B_{z_0}(1)\subset Z, z_0, d_Z), 
(B_{x_0}(1)\subset M_k, x_0, d_{M_k})) < \vare_n=S_n/3.
\end{equation}
So there is a map 
$F_k: B_{x_0}(1) \subset M_k \longmapsto B_{z_0}(1)\subset Z$,
with $F_k(x_0)=z_0$ that is $\vare_n$ almost distance preserving, 
\begin{equation}\label{almodispres}
|d_{M_k}(x_1, x_2)- d_{Z}(F_k(x_1), F_k(x_2))| < \vare_n \qquad
\forall x_1, x_2 \in X,
\end{equation}
and
$\vare_n$ almost onto.  
\begin{equation}\label{almoonto}
\forall z \in B_{z_0}(1)\subset Z,
\,\,\,\, \exists x_z\in B_{x_0}(1)\subset M_k\,\, s.t.\,\,
d_{Z}(F_k(x_z), x)<\vare_n.
\end{equation}
The map $F_k$ can also be thought of as a map defined on 
$B_{x_0}(d_k) \in M^n$.

Note that $F_k$ maps the halfway points, $\gamma_k(d_k/2)$, into
an annulus, 
$$
F_k(\gamma_k(d_k/2)) \subset Ann_{z_0}(1/2 -\vare_n, 1/2 +\vare_n).
$$
Thus there is subsequence of the $F_k(\gamma_k(d_k/2))$ that converges
to a point 
\begin{equation}\label{eqnz1}
z_1 \in Cl(Ann_{z_0}(1/2 -\vare_n, 1/2 +\vare_n)).
\end{equation}
In particular, we can choose a $k \ge N$ such that
\begin{equation}\label{eqnz1k}
d_Z(F_k(\gamma_k(d_k/2)), z_1) < \vare_n.
\end{equation}
We will show that this $z_1$ is our ``cut point''.  To find the cut point
in $Y$, just set $Z=Y$.

We claim that $z_1$ has no ray based at $z_0$ passing through it.  Suppose,
on the contrary, that it does.  Then there is a curve,
$C:[0,1) \longmapsto Z$ such that $c(0)=0$, $c(t_1)=z_1$ and
%\begin{equation} \label{Cmin}
$d_Y(c(t), c(s))=|s-t| \qquad \forall \, t,\,s\,\in [0, 1)$.
%\end{equation}

Let $z_2=C(1/2 + h)$ where $h \in [3 S_n, 1/2)$.  By (\ref{eqnz1}),
$t_1 \ge 1/2-\vare_n$, so
\begin{equation}\label{notpolecut}
d_Z(z_2, z_1) = (1/2 +h) - t_1 \le (1/2+h) - (1/2 -\vare_n)= h+\vare_n. 
\end{equation} 
By (\ref{almoonto}), there exist $x_2 \in B_{x_0}(1)\subset M_k$
which is mapped almost onto $z_2$,
\begin{equation}\label{eqnz2}
d_Z(F_k(x_2), z_2) < \vare_n.
\end{equation}
By (\ref{almodispres}), the triangle inequality and (\ref{eqnz2}), we know that
\begin{eqnarray*}
d_{M}(x_0, x_2)& =& d_{M_k}(x_0, x_2) d_k 
                > (d_{Z}(F_k(x_0), F_k(x_2)) - \vare_n)d_k  \\
&  \ge & (d_{Z}(z_0, z_2) - d_{Z}(z_2, F_k(x_2)) - \vare_n)d_k
    > ((1/2 +h) - 2 \vare_n)d_k.
\end{eqnarray*}
Recall that $h \ge 3S_n$ and $\vare_n =S_n/3$, so 
$x_2 \in \partial B_{x_0}(rd_k)$ with 
$r \ge 1/2 +S_n$.  Thus by the uniform cut property, (\ref{lemmaseqn}),
\begin{equation}\label{lemmaseqn2}
d_M(x_2, \gamma_k(d_k/2))\ge  (r-1/2+2 S_n) d_k
\ge (h-2\vare_n + 2S_n)d_k.
\end{equation}

However, $F_k$ is $\vare_n$ almost distance preserving, (\ref{almodispres}). 
So applying the triangle inequality, (\ref{eqnz2}), (\ref{notpolecut}),
and (\ref{eqnz1k}), we have
\begin{eqnarray*}
d_M(x_2, \gamma_k(d_k/2)) & = & d_{M_k}(x_2, \gamma_k(d_k/2))d_k  \\
&<& (d_{Z}(F_k(x_2), F_k(\gamma_k(d_k/2))) + \vare_n) d_k  \\
& \le& (d_{Z}(F_k(x_2), z_2) + d_Z(z_2, z_1) + d_Z(z_1,F_k(\gamma_k(d_k/2)))
   + \vare_n) d_k  \\
&<& (\vare_n + (h+\vare_n) + \vare_n + \vare_n) d_k.
\end{eqnarray*}
Combining this equation with (\ref{lemmaseqn2}), we get
$$
(h + \vare_n + 3\vare_n) d_k > (h-2\vare_n + 2S_n)d_k.
$$
So $\vare_n > S_n/3$, contradicting (\ref{eqnvareSn}).
\ProofEnd
\end{section}

%%%%%%%%%%%%%%%%%%%%%%%%%%%%%%%%%%%%%%%%%%%%%%%%%%%%%%%%%%%%%%%%%%%%%%%%%%

\end{document}